\documentclass[11pt,a4wide,oneside,reqno]{amsart}
\usepackage{amssymb,amsmath,url,a4wide}
\usepackage[utf8]{inputenc}
\usepackage{header}
\usepackage[left=3cm,right=3cm,top=3cm,bottom=2cm]{geometry}
\usepackage{enumitem}

\newcommand{\conf}{\mathcal{C}}
\newcommand{\niconf}{\mathcal{C^\#}}
\newcommand{\iconf}{\mathcal{C}^\mathrm{inj}}
\newcommand{\g}{\mathfrak{D}}
\newcommand{\citep}[1]{\cite{#1}}

\setlist[enumerate]{label=(\roman*)}
\numberwithin{equation}{section}
\usepackage[url=false,doi=true,maxnames=4]{biblatex}
\usepackage{xpatch}
\xpatchbibmacro{volume+number+eid}{%
  \setunit*{\adddot}%
}{%
}{}{}
\DeclareFieldFormat[article]{number}{\mkbibparens{#1}}
\renewbibmacro{in:}{%
  \ifentrytype{article}{}{\printtext{\bibstring{in}\intitlepunct}}}
\title{Common and Sidorenko equations in Abelian groups}
\author{Leo Versteegen}
\address{Department of Pure Mathematics and Mathematical Statistics, Centre for
Mathematical Sciences, Wilberforce Road, Cambridge CB3 0WB, United Kingdom}
\email{lvv23@dpmms.cam.ac.uk}
\begin{document}
\begin{abstract}
\noindent A linear configuration is said to be common in a finite Abelian group $G$ if for every 2-coloring of $G$ the number of monochromatic instances of the configuration is at least as large as for a randomly chosen coloring. Saad and Wolf conjectured that if a configuration is defined as the solution set of a single homogeneous equation over $G$, then it is common in $\F_p^n$ if and only if the equation's coefficients can be partitioned into pairs that sum to zero mod $p$. This was proven by Fox, Pham and Zhao for sufficiently large $n$. We generalize their result to all sufficiently large Abelian groups $G$ for which the equation's coefficients are coprime to $\lv G\rv$.
\end{abstract}
\maketitle
\section{Introduction}
Let $G$ be a finite Abelian group. For a pair of integers $k\leq d$, we consider matrices $L\in \Z^{k\times d}$ that have full rank in $G$, in the sense that the maps $G^d\rightarrow G^k$ the matrices induce are surjective. The \textit{(linear) configuration given by $L$}, denoted by $\conf(L)$, is the kernel of this map in $G^d$. We refer to the elements of $\conf(L)$ as the \textit{instances} of the configuration. For a function $f\colon G \rightarrow \R$ we define the \textit{arithmetic multiplicity} as\footnote{Here and throughout the paper we write for a set $X$ and a function $f\colon X\rightarrow \C$ the \textit{average of $f$ over $X$} as $\E(f)=\E_{x\in X} f(x)=\frac{1}{\lv X\rv}\sum_{x\in X} f(x)$.}
\begin{align*}
t_L(f)=\E_{v\in \conf(L)} \prod_{i=1}^d f(v_i)=\frac{1}{\vert \conf(L)\vert} \sum_{v\in \conf(L)} \prod_{i=1}^d f(v_i).
\end{align*}
For a set $A\subset G$ we denote by $1_A$ the indicator function of $A$ and write $t_L(A)$ to mean $t_L(1_A)$. In analogy to graph-theoretical concepts, Saad and Wolf \cite{saad-wolf-sidorenko-common} introduced the notions of a configuration being \textit{Sidorenko} or \textit{common} (see the references in the introduction of \cite{saad-wolf-sidorenko-common}). We define commonness in a slightly different way here, which we believe is better suited to general finite Abelian groups instead of the specific classes of groups, which were the primary focus of Saad and Wolf.
\begin{definition}
We call $L$ \textit{fully Sidorenko in $G$} if for all $f\colon G\rightarrow [0,1]$
\begin{align*}
t_L(f)\geq \E(f)^d.
\end{align*}
Furthermore, $L$ is simply called \textit{Sidorenko in $G$} if for all $A\subset G$
\begin{align*}
t_L(A)\geq \lp \frac{\vert A\vert}{\vert G\vert}\rp^d.
\end{align*}
We say that $L$ is \textit{fully common in $G$} if for all $f\colon G\rightarrow [0,1]$ 
\begin{align*}
t_L(f)+t_L(1-f)\geq \lp\frac{1}{2}\rp^{k-1}.
\end{align*}
Finally, $L$ is said to be \textit{common in $G$} if for all $A\subset G$
\begin{align*}
t_{L}(A) + t_{L}(A^C) \geq \lp\frac{1}{2}\rp^{k-1}.
\end{align*}
\end{definition}
Clearly, if a configuration is fully common or fully Sidorenko, then it is also common or Sidorenko, respectively. Conversely, if the group $G$ is sufficiently large and $L$ is non-degenerate, a standard argument allows us to convert functions $f\colon G\rightarrow [0,1]$ into sets $A\subset G$ such that $t_\conf(f)\approx t_\conf(1_A)$.\\\\
As a first step towards a general characterization of Sidorenko and common configurations, Saad and Wolf laid out in \cite{saad-wolf-sidorenko-common} some initial results, conjectures and questions about necessary and sufficient conditions for a configuration to be Sidorenko or common. Among these was the following.
\begin{conjecture}\label{conj:saad-wolf}
Let $d\geq 2$ even. A linear configuration defined by a single equation in $d$ variables, i.e., as the kernel of a matrix $L\in \Z^{1\times d}$, is common in $\F_p^n$ if and only if the set $[d]:=\{1,\ldots,d\}$ can be partitioned into pairs such that for each pair $\{i,j\}$ we have $L_{1i}\equiv -L_{1j}$ $\mod p$.
\end{conjecture}
Saad and Wolf remarked that one direction is easy: from the assumption about the coefficients of $L$, it follows by an application of the Cauchy-Schwarz inequality that the configuration is in fact Sidorenko. In \cite{fox-pham-zhao-sidorenko-equations}, Fox, Pham and Zhao went on to show the reverse direction, thus proving the conjecture. In this paper, we generalize these results from $\F_p^n$ to all finite Abelian groups $G$ in which the coefficients of $L$ are coprime to $\lv G\rv$. First we need to make precise how the condition $L_{1i}\equiv -L_{1j}$ $\mod p$ translates to general finite Abelian groups.
\begin{definition}[Canceling partitions]
A pair of integers $(m,n)$ is \textit{canceling in $G$} if their sum acts trivially on the group, i.e., $(m+n)x=0$ for all $x\in G$. For a matrix $L\in \Z^{k\times d}$, a \textit{canceling partition in $G$} is a partition of $[d]$ into (distinct) pairs such that for each pair $\{i,j\}$ and each $h\in [k]$ the pair $(L_{hi},L_{hj})$ is canceling in $G$.
\end{definition}
The following result settles which single equations are fully common and which are fully Sidorenko under a straightforward condition on the non-degeneracy of the equations.
\begin{restatable}{theorem}{thmSingleEquation}\label{thm:single-equation}
Let $G$ be a finite Abelian group and let $L\in \Z^{1\times d}$ be such that all coefficients in $L$ are coprime to the order of $G$.
\begin{enumerate}
\item If $d$ is even and a canceling partition exists for $L$, then $L$ is fully Sidorenko\footnote{It is easy to see that this holds also for $k\geq 1$.}.
\item If $d$ is odd, then $L$ is fully common in $G$.
\item If $d$ is odd, then $L$ is not fully Sidorenko in $G$.
\item If $d$ is even but no canceling partition exists for $L$, then $L$ is not fully common (and hence not fully Sidorenko) in $G$.
\end{enumerate}
\end{restatable}
Part (i) was observed for $G=\F_p^n$ by Saad and Wolf and their proof carries over directly. The commonness of a single equation in an odd number of variables, part (ii) of the above theorem, was already proved in \cite{odd-variable-equations-common} by Cameron, Cilleruelo and Serra, even for non-Abelian groups, but we include a proof for completeness. The remaining two statements were proved in \cite{fox-pham-zhao-sidorenko-equations} for $G=\F_p^n$. Our proof resembles that in \cite{fox-pham-zhao-sidorenko-equations} in so far as we also construct a suitable function via its Fourier coefficients. However, the possibility of $G$ having no small subgroup means we have to exercise considerably more care when choosing these coefficients, as we will explain in our concluding remarks.\\\\
As a corollary to \Cref{thm:single-equation} we obtain the following result in terms of sets rather than functions.
\begin{restatable}{corollary}{corSingleEquation}\label{cor:single-equation}
For all $d\in \N$, there exists a constant $C$ such that the following is true. Let $G$ be a group with $\vert G\vert >C$, and let $L\in \Z^{1\times d}$ be such that all coefficients in $L$ are coprime to the order of $G$. If $d$ is odd, then $L$ is not Sidorenko in $G$. If $d$ is even but $L$ has no canceling partition in $G$, then $L$ is not common in $G$.
\end{restatable}
Conceivably, the restriction that the coefficients should be coprime to the order of the group could be weakened in some situations. For the groups $\Z_p$ or $\F_p^n$, this would mean that we allow (some) coefficients of the equation to be zero. This adds nothing of interest because such an equation will be Sidorenko or common if and only if the reduced equation (with all zero coefficients removed) is. For other finite Abelian groups, the problem is more delicate, but any formulation of a result seems to require a large of number of case distinctions regarding the group decomposition, rendering it of limited interest.
\subsection*{Acknowledgments}
The author is grateful to be funded by Trinity College of the University of Cambridge through the Trinity External Researcher Studentship. Further, he wishes to express his deep gratitude to his master thesis supervisor Mathias Schacht  for introducing him to the area of additive combinatorics, to Emil Powierski for many hours of enlightening discussion on the topic and to his PhD supervisor  Julia Wolf for her generous help in improving this paper.
\section{The proof of \Cref{thm:single-equation}}
Our proof of \Cref{thm:single-equation} relies heavily on the Fourier transform. For the convenience of the reader, we recall briefly how it is defined in a discrete setting and state without proof its most important properties.\\\\
For a finite Abelian group $G$, consider the set $\hat{G}$ of \textit{characters}, i.e., group homomorphisms $\gamma \colon G\rightarrow \{z\in \C:\lv z\rv=1\}$. Equipped with pointwise multiplication as a binary operation between characters, $\hat{G}$ becomes a group which is is isomorphic to $G$. For later use, we fix an arbitrary isomorphism $\g\colon G\rightarrow \hat{G}$. One reason characters are so useful to us is that they detect the additive identity of a group by the following well-known identity.
\begin{align}\label{eq:char-ortho-sum}
\sum_{\gamma\in \hat{G}} \gamma(x)=\begin{cases}\vert G\vert &\text{if } x=0,\\0&\text{otherwise.}\end{cases}
\end{align}
The other benefit of characters is that we may use them to analyze functions via their Fourier coefficients. Recall that for a function $f\colon G\rightarrow \C$ the \textit{Fourier transform of $f$} is the function
\begin{align*}
\hat{f}\colon \hat{G} \rightarrow \C \qquad \gamma \mapsto \hat{f}(\gamma)=\E_{x\in G} \gamma(x)f(x).
\end{align*}
The operator that maps a function to its Fourier transform is linear and a scaled isometry in the sense that for all $f_1,f_2\in \Map(G,\C)$ it holds that
\begin{align*}
\sum_{\gamma\in \hat{G}} \hat{f_1}(\gamma)\ov{\hat{f_2}(\gamma)}=\E_{x\in G} f_1(x)\ov{f_2(x)}.
\end{align*}
Furthermore, we have for all $f$ the inversion formula 
\begin{align*}
f(x)=\sum_{\gamma\in \hat{G}} \hat{f}(\gamma)\gamma(-x).
\end{align*}
Real-valued functions $f$ do not generally have a real-valued Fourier transform. If and only $f$ is real, then $\hat{f}(\gamma\inv)=\ov{\hat{f}(\gamma)}$ for every $\gamma\in \hat{G}$.\\\\
The identity element of $\hat{G}$ is the character that always takes the value 1 and it plays a special role with respect to the Fourier transform. Namely, for every function $f\colon G\rightarrow \C$ we have that $\hat{f}(1)=\E(f)$.\\\\
The Fourier transform allows us to express the number of solutions to a system of equations in a concise manner.
\begin{lemma}\label{lemma:fourier-density}
Let $L\in \Z^{1\times d}$ have full rank in $G$. Then we have for all $f\colon G\rightarrow [0,1]$
\begin{align}\label{eq:full-rank-hom-density}
t_{L}(f)=\sum_{\gamma\in \hat{G}}\prod_{i=1}^d \hat{f}(\gamma^{L_{1i}})=\sum_{x\in G}\prod_{i=1}^d \hat{f}\lp \g\lp L_{1i}x\rp\rp.
\end{align}
\end{lemma}
\begin{proof}
Equation \eqref{eq:char-ortho-sum} allows us to express the fact that a vector $v\in G^d$ is a solution of $Lv=0$ via
\begin{align*}
\sum_{\gamma\in \hat{G}} \gamma(Lv)=\begin{cases}\vert G\vert &\text{if } Lv=0,\\ 0&\text{otherwise.}\end{cases}
\end{align*}
Using that $\lv\{v\in G^d:Lv=0\}\rv = \lv G\rv^{d-1}$, this means we can rewrite the arithmetic multiplicity as
\begin{align*}
\E_{\substack{v\in G^d\\Lv=0}} \prod_{i=1}^d f(v_i)&=\E_{v\in G^d}\lp \prod_{i=1}^d f(v_i) \rp\lp\sum_{\gamma\in \hat{G}} \gamma(Lv) \rp=\sum_{\gamma\in \hat{G}}\E_{v\in G^d} \prod_{i=1}^d f(v_i) \gamma^{L_{1i}}( v_i).
\end{align*}
If we split the expectation into the different coordinates, we see that by definition of the Fourier transform this is the same as
\begin{align*}
\sum_{\gamma\in \hat{G}}\prod_{i=1}^d \hat{f}(\gamma^{L_{1i}}),
\end{align*}
as desired.
\end{proof}
We are now ready to prove \Cref{thm:single-equation}, which we restate here for the convenience of the reader.
\thmSingleEquation*
\begin{proof}
It is easy to see that systems consisting of only one equation have full rank if at least one coefficient of $L$ is coprime to $\vert G\vert$. This is the case here, so $L$ has full rank and we may apply \Cref{lemma:fourier-density} to see
\begin{align}\label{eq:single-density}
t_{L}(f)=(\E(f))^d+\sum_{v\in G\setminus \{0\}}\prod_{i=1}^d \hat{f}\lp \g(L_{1i}x)\rp.
\end{align}
Throughout the proof, we shall refer to the sum over the non-zero values of $v$ on the right-hand side as the \textit{deviation} for $f$. If a canceling partition exists, then the deviation is always positive because $\hat{f}\lp \g(L_{1i}x)\rp\hat{f}\lp \g(-L_{1i}x)\rp=\vert\hat{f}\lp \g(L_{1i}x)\rp\vert^2$, which proves (i). On the other hand, proving that $L$ is not fully Sidorenko is equivalent to finding $f$ such that the deviation is negative. Before constructing such a function $f$, we will ascertain what is needed to prove or disprove commonness.\\\\
By linearity of the Fourier transform, for any $f\colon G\rightarrow \C$ we have $\hat{1-f}=\hat{1}-\hat{f}$. It is easy to see by \eqref{eq:char-ortho-sum} that for $x\neq 0$, we have $\hat{1}(\g(x))=0$. Hence $\hat{(1-f)}(\g(x))=-\hat{f}(\g(x))$. Also, since all $L_{1i}$ are coprime to the order of the group, we can argue that $L_{1i}x$ can only be zero if $x$ is zero. Hence, $\hat{(1-f)}(\g(L_{1i}x))=-\hat{f}(\g(L_{1i}x))$ for $x\neq 0$. This gives
\begin{align}\label{eq:double-density}
t_{L}(f)+t_{L}(1-f)&=(\E(f))^d + (1-\E(f))^d\nonumber\\&\hspace*{2.0cm}+\sum_{x\neq 0}\lp\prod_{i=1}^d \hat{f}\lp \g(L_{1i}x)\rp+\prod_{i=1}^d \lp-\hat{f}\lp \g(L_{1i}x)\rp\rp\rp.
\end{align} 
Clearly, if $d$ is odd, then the two products in the sum on the right cancel. The system is thereby shown to be common as the term $(\E(f))^d+(1-\E(f))^d$ is minimized by $\E(f)=\frac{1}{2}$. This proves (ii). If $d$ is even, \eqref{eq:double-density} becomes
\begin{align}\label{eq:double-density-cpct}
t_{L}(f)+t_{L}(1-f)&=(\E(f))^d + (1-\E(f))^d+2\sum_{x\neq 0}\prod_{i=1}^d \hat{f}\lp \g(L_{1i}x)\rp.
\end{align}
The sum on the right is of course again the deviation from \eqref{eq:single-density}. Thus, to prove both our remaining claims it suffices to construct, for a given $L$ that lacks a canceling partition (which is always the case when $d$ is odd), a function $f\colon G\rightarrow [0,1]$ with $\E(f)=\frac{1}{2}$ such that the deviation is negative. We will construct $f$ as the Fourier-inverse of an appropriate function.\\\\
To this end, pick an element $a\in G$ of maximal order. Observe that the order of any other element in $G$ divides that of $a$. Note that the order of $a$ cannot be 2, because if all elements had order two then every partition of the coefficients of $L$ into pairs would be cancelling. The order of $a$ is also coprime to every coefficient of $L$. Therefore, no element of $U= \{L_{11}a,\ldots,L_{1d}a\}$ is its own inverse, in particular $0\notin U$. Let now $S=\{s_1,\ldots,s_r\} \subset U$ be such that for each $b\in U$ exactly one of $b$ and $-b$ is in $S$. For each $\vp\in \R$, define a map
\begin{align*}
g_\vp \colon G\rightarrow \C \qquad y \mapsto 
\begin{cases}
\frac{1}{4r}\ex\lp \frac{\vp}{(2d)^j} \rp&\text{if } y=s_j,\\
\frac{1}{4r}\ex\lp -\frac{\vp}{(2d)^j}\rp&\text{if } y=-s_j,\\
\frac{1}{2} &\text{if } y=0,\\
0&\text{otherwise,}
\end{cases}
\end{align*}
where we have written $\ex(\cdot)$ for $\exp(2\pi i \cdot)$. This map is well defined by definition of $S$ and has support $U\cup \{0\}$. We now define $f_\vp$ to be the Fourier inverse of $g_\vp$. To be precise, we define
\begin{align*}
f_\vp\colon G\rightarrow \C \qquad x \mapsto \sum_{\gamma\in \hat{G}} g_\vp(\g\inv(\gamma)) \gamma(-x).
\end{align*}
Because $g_\vp(\g\inv(\gamma\inv))=\ov{g_\vp(\g\inv(\gamma))}$ for all $\gamma\in \hat{G}$, $f_\vp$ takes only real values. Moreover, the image of $f_\vp$ is contained in $[0,1]$. Indeed, for~$x\in G$,
\begin{align*}
\lv f_\vp(x)-\frac{1}{2}\rv&=\lv \sum_{\gamma\neq 1} g_\vp(\g\inv(\gamma))\gamma(-x)\rv\leq \sum_{i=1}^r \vert g_\vp(s_i)\vert + \vert g_\vp(-s_i)\vert=\frac{1}{2}.
\end{align*}
Furthermore, the average of $f_\vp$ is $g_\vp(0)=\frac{1}{2}$ and $\hat{f}=g_\vp \circ \g\inv$. We now wish to find $\vp$ such that the deviation for $f_\vp$ is negative.\\\\
Observe first that if a summand $\prod_{i=1}^d \hat{f_\vp}\lp \g(L_{1i}v)\rp$ is non-zero, it will be $\lp \frac{1}{4r}\rp^d$ times a complex number on the unit circle. Let $X\subset G\setminus\{0\}$ be the set of $v$ for which this is the case. Note that~$X$ is non-empty, as $g_\vp$ is designed so that at least $a\in X$. Now, fix $x\in X$. We have $L_{11}x=\pm L_{1i}a$ for some $i\in [d]$. Therefore, if $mx=0$ for some $m\in \Z$, we have $mL_{1i}a=0$. But since $L_{1i}$ is coprime to the order of $a$, we get that $ma=0$. Hence $x$ is also of maximal order in $G$ and the order of every element in $G$ divides that of $x$.\\\\
Each factor $\hat{f}(L_{1i}x)$ is non-zero, so that for each $i\in [k]$ there are $j_i\in [r]$ and $\sigma_i\in \{\pm 1\}$ such that $L_{1i}x=\sigma_i \cdot s_{j_i}$. The summand $\prod_{i=1}^d \hat{f_\vp}\lp \g(L_{1i}x)\rp$ can therefore be written as
\begin{align*}
\prod_{i=1}^d \hat{f_\vp}\lp \g(L_{1i}x)\rp=\prod_{i=1}^d g_\vp\lp L_{1i}x\rp=\lp \frac{1}{4r}\rp^d \ex \lp \vp \sum_{i=1}^d \sigma_i \frac{1}{(2d)^{j_i}}\rp=\lp \frac{1}{4r}\rp^d \ex \lp \vp \sum_{j=1}^r \sum_{\substack{i\in [d]\\j_i=j}} \sigma_i \frac{1}{(2d)^j}\rp.
\end{align*}
We claim that there is an index $j$ such that
\begin{align*}
\lv \{i\in [d]: j_i=j, \sigma_i=1\}\rv \neq \lv\{i\in [d]: j_i=j, \sigma_i=-1\}\rv,
\end{align*}
for otherwise we could partition $[d]$ into pairs such that for every pair $\{i_1,i_2\}$ among them $j_{i_1}=j_{i_2}$ and $\sigma_{i_1}=-\sigma_{i_2}$, implying that $(L_{1i_1}+L_{1i_2})x=0$. But because the order of every element in $G$ divides the order of $x$, we would have $(L_{1i_1}+L_{1i_2})y=0$ for all $y\in G$, i.e., there would be a partition of $[d]$ into canceling pairs in $G$, which we have assumed not to be the case.\\\\
Choose the minimal index $j$ with this property and denote it by $j_*$. We have
\begin{align*}
\lv \sum_{j=1}^r \sum_{\substack{i\in [d]\\j_i=j}} \sigma_i \frac{1}{(2d)^j}\rv &\geq \lv \sum_{\substack{i\in [d]\\j_i=j_*}} \sigma_i\frac{1}{(2d)^{j_*}}\rv-\sum_{j=j_*+1}^r \sum_{\substack{i\in [d]\\j_i=j}} \frac{1}{(2d)^j}\\
&\geq \frac{1}{(2d)^{j_*}}-\frac{1}{2^{j_*+1}d^{j_*}}=\frac{1}{2^{j_*+1}d^{j_*}}>0.
\end{align*}
Let
\begin{align*}
c_x:=\sum_{i=1}^d \sigma_i \frac{1}{(2d)^{j_i}}= \sum_{j=1}^r \sum_{\substack{i\in [d]\\j_i=j}} \sigma_i \frac{1}{(2d)^j}.
\end{align*}
Since for each $j\in [d]$ the summand on the right-hand side is an integer multiple of $(2d)^{-r}$, $c_x$ will be a multiple of $(2d)^{-r}$ also. Because $c_x\neq 0$, $c_x$ must have modulus at least $(2d)^{-r}$. On the other hand, the triangle inequality gives $\vert c_x\vert \leq \frac{1}{2}$.\\\\
Consider now the function $\psi$ that maps a phase $\vp$ to $(4r)^d$ times the deviation of $f_\vp$, i.e.,
\begin{align*}
\psi \colon [0,(2d)^r] \rightarrow \C \qquad \vp \mapsto (4r)^d\sum_{x\neq 0}\prod_{i=1}^d g_\vp\lp L_{1i}x\rp=\sum_{x\in X} \ex\lp c_x\vp\rp.
\end{align*}
Note that $\vert \psi\vert$ is bounded from above by $\vert X\vert$. If $\vp$ is sufficiently close to the boundaries of the interval $[0,(2d)^r]$, $\psi(\vp)$ will have positive real part. To be precise, if $\vp\leq\frac{1}{4}$ or $\vp \geq (2d)^r-\frac{1}{4}$, then
\begin{align*}
\Re\lp\psi(\vp)\rp= \sum_{x\in X} \cos(2\pi c_x \vp) \geq \frac{\vert X\vert}{\sqrt{2}}.
\end{align*}
On the other hand, the average of $\psi$ is 0 in the sense that
\begin{align*}
\int_0^{(2d)^r} \psi(\vp)\d \vp=\frac{1}{2\pi i}\sum_{x\in X} \frac{1}{c_x}\left[\ex(c_x\vp)\right]_0^{(2d)^r}=0.
\end{align*}
Now let 
\begin{align*}
A&=\left\{\vp \in [0,(2d)^r]:\Re(\psi(\vp))\leq - \frac{\vert X\vert}{2\sqrt{2}(2d)^r}\right\},\\
B&=\left\{\vp \in [0,(2d)^r]:- \frac{\vert X\vert}{2\sqrt{2}(2d)^r}<\Re(\psi(\vp))\leq 0\right\}.
\end{align*} 
Denoting by $\lambda$ the Lebesgue measure on $[0,(2d)^r]$ we obtain
\begin{align*}
\Re\lp\int_0^{(2d)^r} \psi(\vp)\d \vp\rp &\geq \frac{1}{2}\cdot \frac{\vert X\vert}{\sqrt{2}} - \lambda(A) \cdot \vert X\vert - \lambda(B)\cdot \frac{\vert X\vert}{2\sqrt{2}(2d)^r}\\
&> \vert X\vert \lp \frac{1}{2\sqrt{2}} - \lambda(A) - \frac{(2d)^r}{2\sqrt{2}(2d)^r}\rp\\
&=-\lambda(A)\cdot \vert X\vert.
\end{align*}
Because we had already established that the integral on the left-hand side is 0, $A$ must have positive measure, meaning that $A$ is not empty. But by the definitions of $\psi$ and $A$, $f_\vp$ has negative deviation for any $\vp\in A$, completing the proof.
\end{proof}
Note that the proof in fact yields an explicit upper bound on the minimal deviation in \eqref{eq:single-density} over all $f\colon G\rightarrow [0,1]$ with $\E(f)=1/2$. Namely, using that $r\leq d$, we have for at least one value of $\vp$
\begin{align}\label{eq:deviation}
\sum_{x\neq 0}\prod_{i=1}^d \hat{f_\vp}\lp \g(L_{1i}x)\rp\leq -\frac{1}{(4r)^d}\cdot\frac{1}{2\sqrt{2}(2d)^r}\leq -\frac{1}{2^{3d+1}\sqrt{2}d^{2d}}.
\end{align}
\section{Deriving \Cref{cor:single-equation} from \Cref{thm:single-equation}}
There are different ways to deduce \Cref{cor:single-equation} from \Cref{thm:single-equation}, but all of them will require that the set 
\begin{align*}
\niconf(L)=\{v\in \conf(L)\subset G^d: \exists i,j\in [d]\colon i\neq j \land v_i=v_j\}
\end{align*}
of \textit{non-injective instances} is not too large. The following lemma, which is inspired by Lemma 2.1 in \cite{reiher-odd-cycles-locally-dense}, makes the dependence on the size of $\niconf(L)$ very explicit.
\begin{lemma}[Conversion of maps to sets]\label{lem:maps-to-sets}
Let $G$ be an Abelian group, $L\in \Z^{k\times d}$ a system of equations and $f\colon G\rightarrow [0,1]$.
\begin{enumerate}
\item There exists a set $A\subset G$ with $\vert A\vert \geq \E(f)\cdot \vert G\vert-1$ and
\begin{align*}
t_L(A) \leq t_L(f) +  \frac{\vert \niconf(L) \cap A^d\vert}{\vert \conf(L)\vert}.
\end{align*}
\item There exists a set $A\subset G$ with
\begin{align*}
t_L(A) + t_L(A^C) \leq t_L(f) + t_L(1-f) + \frac{\vert \niconf(L)\cap A^d\vert+\vert \niconf(L) \cap (A^C)^d\vert}{\vert\conf(L) \vert}.
\end{align*}
\end{enumerate}
\end{lemma}
\begin{proof}
We prove only (i) since (ii) can be shown by an argument that is almost identical. Let $\iconf=\conf(L)\setminus \niconf(L)$. The map
\begin{align*}
\psi \colon \Map(G,[0,1]) \rightarrow \R \qquad g\mapsto \sum_{v\in \iconf} \prod_{i=1}^k g(v_i)
\end{align*}
attains its minimum on the compact and non-empty set $\{ g\in [0,1]^G:\E(g)=\E(f)\}$. Of all the maps $g$ achieving this minimum, we denote by $g_0$ one for which also the set $R=\{x\in G:g(x)\notin \{0,1\}\}$ is smallest. We show that $R$ can contain at most one element. Assume to the contrary that there are two different elements $a,b\in R$. For $\eta\in \R$, consider
\begin{align*}
g_\eta\colon G\rightarrow \R \qquad x\mapsto \begin{cases}g_0(a)+\eta &\text{if } x=a,\\g_0(b)-\eta &\text{if } x=b,\\g_0(x) &\text{otherwise.}\end{cases}
\end{align*}
Clearly, $\E(g_\eta)=\E(f)$ and because $a,b\in R$ we have $\im (g_\eta)\subset [0,1]$ for small enough $\eta$. It is easy to see that $\eta\mapsto \psi(g_\eta)$ is a quadratic polynomial with constant term $\psi(g_0)$ and a non-positive quadratic coefficient. By minimality of $\psi(g_0)$, this function must then be constant. But then we may choose $\eta$ so that at least one of $a$ and $b$ takes a value in $\{0,1\}$ under $g_\eta$. This contradicts the minimality of the set $R$.\\\\
Now let $A=g_0\inv(1)$. Then $1_A(x)=g_0(x)$ for $x\notin R$ and because $\vert R\vert \leq 1$, we have that $\vert A\vert>\E(f)\cdot\vert G\vert-1$. On the other hand, the monotonicity of $\psi$ gives $\psi(1_A)\leq \psi(g_0)$, and $\psi(g_0)$ in turn is less than $\psi(f)$ by the minimality of $g_0$ under $\psi$ among the functions with average $\E(f)$. All in all, we get
\begin{align*}
t_L(A)&=\frac{1}{\vert \conf(L)\vert}\lp \psi(1_A) + \sum_{v\in \niconf(L)} \prod_{i=1}^k 1_A(v_i) \rp\leq \frac{\psi(f) +\vert \niconf \cap A^d\vert}{\lv\conf(L)\rv} \leq t_L(f) +  \frac{\vert \niconf \cap A^d\vert}{\vert \conf(L)\vert}.
\end{align*}
\end{proof}
Unfortunately, there exist configurations in which all instances are non-injective. Having said that, the following lemma provides a useful bound on the size of $\niconf(L)$ for all matrices $L$ consisting of a single row under a rather weak assumption about the coefficients of the equation.
\begin{lemma}\label{lemma:non-injective-single-equation}
Let $G$ be any finite Abelian group, $d\geq 3$ and $L\in \Z^{1\times d}$. If at least three coefficients of $L$ are coprime to $\vert G\vert$, then
\begin{align*}
\vert\niconf(L)\vert \leq \binom{d}{2}\frac{\vert\conf(L)\vert}{\vert G\vert}.
\end{align*}
\end{lemma}
\begin{proof}
The matrix $L$ describes a surjective linear map. Indeed, if $L_{1i}$ is coprime to $\lv G\rv$ for some $i\in [d]$, then fixing all but the $i$th coordinate of $v\mapsto Lv$ defines an automorphism of $G$. This means $\vert \conf(L)\vert=\lv \ker L\rv=\vert G\vert^{d-1}$. It is now enough to show that for any pair $\{i,j\}\in[d]^{(2)}$ the number of solutions where $v_i=v_j$ can be bounded by $\vert G\vert^{d-2}$. Indeed, for each such pair there is, by assumption, a third index $m$ such that~$L_{1m}$ is coprime to $\vert G\vert$. Then, if we choose all coordinates except $v_j$ and $v_m$ freely from $G$ and set $v_j=v_i$, there is exactly one $v_m\in G$ such that $Lv=0$.
\end{proof}
Since in \Cref{thm:single-equation} we assume that all coefficients are coprime to the order of the group, we may apply \Cref{lemma:non-injective-single-equation} to deduce \Cref{cor:single-equation}. We recall its statement here.
\corSingleEquation*
\begin{proof}
We denote by $\Delta$ the modulus of the deviation achieved in \eqref{eq:deviation}, i.e.,
\begin{align*}
\Delta = \frac{1}{2^{3d+1}\sqrt{2}d^{2d}}.
\end{align*}
Choose $C= d^2/\Delta$ and let $G$ and $L$ be as in the statement of \Cref{cor:single-equation}. The proof of \Cref{thm:single-equation} yields a function $f$ with~$\E(f)=1/2$, $t_L(f)\leq 1/2^d-\Delta$ and $t_L(1-f)\leq 1/2^d-\Delta$. Now apply \Cref{lem:maps-to-sets} to obtain sets $A_1,A_2$ such that
\begin{align*}
\vert A_1\vert &\geq \E(f) \vert G\vert -1\\
t_{L}(A_1) &\leq t_{L}(f) +  \frac{\vert \niconf(L)\vert}{\vert \conf(L)\vert}\\
t_{L}(A_2) + t_{L}(A_2^C) &\leq t_{L}(f) + t_{L}(1-f) + \frac{2\vert \niconf(L) \vert}{\vert \conf(L)\vert}.
\end{align*}
Recall \Cref{lemma:non-injective-single-equation}, which stated that matrices $L\in \Z^{1\times d}$ with at least 3 coefficients that are coprime to $\vert G\vert$ satisfy
\begin{align*}
\vert\niconf(L)\vert \leq \binom{d}{2}\frac{\vert\conf(L)\vert}{\vert G\vert}.
\end{align*}
Since a matrix without a canceling partition must have at least three coefficients and the coefficients are coprime to $\lv G\rv$ by assumption, the lemma applies. We have
\begin{align*}
t_{L}(A_2) + t_{L}(A_2^C) &\leq \lp \frac{1}{2}\rp^{d-1} - 2\Delta + \frac{d^2}{C}\leq \lp\frac{1}{2}\rp^{d-1}-\Delta,
\end{align*}
showing that $L$ is not common. In order to show that for odd $d$, $L$ is not Sidorenko in $G$, we need to deal with the fact that $\vert A_1\vert/\vert G\vert$ could be slightly less than $1/2$. We apply the mean value theorem to the monomial $x^d$ to obtain the lower bound
\begin{align*}
\lp\frac{\vert A_1\vert}{\vert G\vert}\rp^d\geq \lp\frac{1}{2}-\frac{1}{\vert G\vert}\rp^d\geq \lp \frac{1}{2}\rp^d - \frac{d}{\vert G\vert}\geq \lp \frac{1}{2}\rp^d - \frac{d}{C},
\end{align*}
whence
\begin{align*}
t_{L}(A_1) &\leq \lp\frac{1}{2}\rp^d- \Delta +  \frac{d^2}{2C} \leq \lp\frac{\vert A_1\vert}{\vert G\vert}\rp^d - \Delta + \frac{d^2+2d}{2C}\leq\lp\frac{\vert A_1\vert}{\vert G\vert}\rp^d-\Delta \lp \frac{d-2}{d} \rp.
\end{align*}
\end{proof}
\section{Concluding remarks}
We conclude the article with a brief comparison between our proof and that given in \cite{fox-pham-zhao-sidorenko-equations} for the special case $G=\F_p^n$, which can easily be reduced to the case $G=\F_p$. There, the witness function~$f$ is also constructed via an inverse Fourier transform. But instead of taking just a few Fourier coefficients to be non-zero, the authors sample for $n=1$ every $\hat{f}(\gamma)$ identically and (essentially) independently for $\gamma\neq 1$ and precompose with a projection to the first coordinate for larger $n$. To ensure that the resulting function $f$ takes values in $[0,1]$, the values $\hat{f}(\gamma)$ have to be scaled appropriately; just as we had to scale by $(4r)\inv$. The key difference is that their scaling factor depends on~$p$, which leads to a much smaller deviation for large $p$. Although we still arrive at a valid proof of \Cref{thm:single-equation}, there are two problems with a deviation that approaches zero as the order of the group grows.\\\\
Firstly, the conversion of $f$ into a set we require to deduce \Cref{cor:single-equation} involves an error term that might move us back above the threshold $\E(f)^d$ if the deviation is too small. If the deviation does not decrease too quickly, a careful conversion might still be possible, but the bound from the proof in \cite{fox-pham-zhao-sidorenko-equations} in combination with \Cref{lem:maps-to-sets} would be insufficient.\\\\
Secondly, and perhaps more importantly, a vanishing deviation would not be in the spirit of \Cref{conj:saad-wolf}. We have mentioned earlier that our definition of commonness differs slightly from that of Saad and Wolf \cite{saad-wolf-sidorenko-common}. By their definition, an equation $L\in \Z^{1\times d}$ is common in $\Z_p$ if
\begin{align*}
\liminf_{p\rightarrow \infty} \min_{A\subset \Z_p} t_L(A)+t_L(A^C) \geq \lp \frac{1}{2}\rp^{d-1}.
\end{align*}
It is therefore a slightly weaker property to be common in $\Z_p$ in the sense of \cite{saad-wolf-sidorenko-common} than to be common in $\Z_p$ for all sufficiently large $p$ in our sense. Conversely, it is a stronger property to be uncommon in the sense of Saad and Wolf. By giving a lower bound on the deviation independent of $p$ when $d$ is even and $L$ has no canceling partition, we have proved this stronger property.
\printbibliography
\end{document}